\documentclass[12pt]{amsart}
\usepackage{amssymb}
\usepackage{amsfonts}
\usepackage{latexsym}
\usepackage{amscd}

\newcommand{\N}{{\mathbb{N}}}

\newcommand{\fl}{\rightarrow}

\newtheorem{lemma}{Lemma}[section]
\newtheorem{corollary}[lemma]{Corollary}
\newtheorem{theorem}[lemma]{Theorem}
\newtheorem{proposition}[lemma]{Proposition}
\newtheorem{remark}[lemma]{Remark}
\newtheorem{definition}[lemma]{Definition}

\begin{document}

\title[Non-unital purely infinite simple rings]{Structure of non-unital purely infinite simple rings.}%
\author{M.A. Gonz\'alez-Barroso}\author{E. Pardo}
\address{Departamento de Matem\'aticas, Universidad de C\'adiz,
Apartado 40, 11510 Puerto Real (C\'adiz), Spain.}
\email{mariangeles.gonzalezbarroso@alum.uca.es}\email{enrique.pardo@uca.es}
\urladdr{http://www.uca.es/dept/matematicas/PPersonales/PardoEspino/EMAIN.HTML}

\thanks{The first author was supported by an FPU fellowship of the Junta de
Andaluc\'{\i}a. Both authors are partially supported by the DGI
and European Regional Development Fund, jointly, through Project
MTM2004-00149, and by PAI III grant FQM-298 of the Junta de
Andaluc\'{\i}a. Also, the second author is partially supported by
the Comissionat per Universitats i Recerca de la Generalitat de
Catalunya.}

\subjclass[2000]{Primary 16D70}

\keywords{Purely infinite simple ring, Zhang's Dichotomy}


\begin{abstract}
In this note, we study the notion of purely infinite simple ring
in the case of non-unital rings, and we obtain an analog to
Zhang's Dichotomy for $\sigma$-unital purely infinite simple
C*-algebras in the purely algebraic context.
\end{abstract}
\maketitle

\section*{Introduction}

In 1981, Cuntz \cite{C2} introduced the concept of a purely
infinite simple C*-algebra. This notion has played a central role
in the development of the theory of C*-algebras in the last two
decades. A large series of contributions, due to Blackadar, Brown,
Lin, Pedersen, Phillips, R\o rdam and Zhang, among others, reflect
the interest in the structure of such algebras. A particular
interest deserves Zhang's result \cite{Zhang}, dividing
$\sigma$-unital purely infinite simple C*-algebras in two types:
unital and stable. This result, known as Zhang's Dichotomy for
$\sigma$-unital purely infinite simple C*-algebras, played a
central role in the study of the structure of corona and
multiplier algebras for C*-algebras with real rank zero.

In 2002, Ara, Goodearl and Pardo \cite{agp} introduced a suitable
definition of a purely infinite simple ring for unital rings,
which agrees with that of Cuntz in the case of C*-algebras, and
studied $K_0$ and $K_1$ groups of a purely infinite simple ring,
specially in the case of von Neumann regular rings lying in this
class. The natural generalization of this definition to the
context of non-unital rings was already considered in \cite{AP},
and also in \cite{exchange}, where Ara showed that every (non
necessarily unital) purely infinite simple ring is an exchange
ring.

In this note, we study the notion of non-unital purely infinite
simple ring considered in \cite{exchange}. We start by comparing
this notion with a different one, inspired in \cite[Theorem
1.6]{agp}, which turns out to be equivalent to the former one for
C*-algebras \cite{C2}, \cite{LZ}. We conclude that the original
definition is stronger that the new one, but it is not clear
whether both definitions are equivalent in the algebraic context.
Finally, using the definition introduced in \cite{exchange}, we
are able to prove an algebraic version of Zhang's result, dividing
$\sigma$-unital purely infinite simple rings in unital and
stable.\vspace{.2truecm}

We need to fix some definitions. Given a ring $R$, we denote by
$M_{\infty }(R)=\varinjlim M_n(R)$, under the maps
$M_n(R)\rightarrow M_{n+1}(R)$ defined by $x\mapsto
\mbox{diag}(x,0)$. Notice that $M_{\infty }(R)$ can also be
described as the ring of countable infinite matrices over $R$ with
only finitely many nonzero entries. Given $p,q\in M_{\infty }(R)$
idempotents, we say that $p$ and $q$ are equivalent, denoted
$p\sim q$, if there exist elements $x,y\in M_{\infty }(R)$ such
that $xy=p$ and $yx=q$. We also write $p\leq q$ provided that
$p=pq=qp$, $p\lesssim q$ if there exists an idempotent $r\in
M_{\infty }(R)$ such that $p\sim r\leq q$, and $p\prec q$ if there
exists an idempotent $r\in M_{\infty }(R)$ such that $p\sim r <q$.
Given idempotents $p,q\in M_{\infty
}(R)$, we define the direct sum of $p$ and $q$ as $p\oplus q=\left(%
\begin{array}{cc}
  p & 0 \\
  0 & q \\
\end{array}%
\right)$. Also, for an idempotent $p\in M_{\infty }(R)$ and a
positive integer $n$, we denote by $n\cdot p$ the direct sum of
$n$ copies of $p$. Two idempotents $e,f$ are said to be
orthogonal, (denoted $e \perp f$) provided that $ef=fe=0$. In that
case, $e+f$ is an idempotent, and $(e+f)R= eR \oplus fR$. An
idempotent $e$ in a ring $R$ is infinite if there exist orthogonal
idempotents $f,g\in R$ such that $e= f+g$ while $e\sim f$ and
$g\ne 0$.

\section{Basic concepts}

In this section we study the  notion of purely infinite simple
ring in the case of non-unital rings. By analogy with the
C*-algebra case, we consider two notions, that turn out to be
equivalent for C*-algebras. The first one is that introduced in
\cite{agp} as a basic definition, and used in \cite{exchange}.

\begin{definition}\label{nounital}
{\rm (\cite[Definition 1.2]{agp}) A ring $R$ is said to be purely
infinite simple if it is simple and every nonzero right ideal
contains an infinite idempotent.}
\end{definition}

The second one is the alternative definition of purely infinite
simple unital ring that rises from \cite[Theorem 1.6]{agp},
adapted to the non-unital case. We borrow the name from \cite[pp.
241--242]{Cohn}.

\begin{definition}\label{weaky}
{\rm A nonzero ring $R$ is $1$-simple if for every nonzero
elements $x,y\in R$ there exist $z,t\in R$ such that $zxt=y$.}
\end{definition}

\begin{remark}\label{observacio 1}

{\rm \noindent {\rm {(1)}} It is easy to see that the definition
of purely infinite simple ring is right-left symmetric.

\noindent {\rm {(2)}} It is clear that, by definition, any
$1$-simple ring is simple.

\noindent {\rm {(3)}} If $R$ is a unital $1$-simple ring, then it
is either a division ring or a purely infinite simple ring
\cite[Theorem 1.6]{agp}.}
\end{remark}

Now, we study the relation between these definitions in the purely
algebraic context.

\begin{proposition}
\label{relacio 1} If $R$ is a (non-unital) purely infinite simple
ring, then it is $1$-simple.
\end{proposition}
\begin{proof}
Let $x,y\in R$ be nonzero elements. By hypothesis, there exists an
infinite idempotent $e\in xR$, so that $e=xr$ for some $r\in R$.
Since $R$ is simple, every nonzero finitely generated projective
module is a generator of the category Mod-$R$. Since $e$ is
infinite and $R$ is simple, it is easy to show that, for any
natural number $n$, there exists a module epimorphism $\varphi _n:
eR\fl n(eR)$. Now, by simplicity, $y\in ReR$, so that
$y=\sum_{i=1}^{m}{z_iet_i}$ for some $z_1, \dots, z_m, t_1, \dots
,t_m\in R$. Hence, multiplication by $(z_1, \dots ,z_m)$ defines a
module homomorphism $\pi :m(eR)\fl R$ such that $y\in
\mbox{Im}(\pi)$. Thus, $\rho =\pi \circ \varphi _m$ defines a
module homomorphism from $eR$ to $R$ such that $y\in
\mbox{Im}(\rho)$. In particular, $y=\rho (et)$ for some $t\in R$.
Since $e=e^2$, for any $a\in R$ we have $\rho (ea)=\rho (e)ea$.
Hence
$$y=\rho(et)=\rho (e)et=\rho (e)x (rt),$$
as desired.
\end{proof}

The converse of Proposition \ref{relacio 1} holds whenever $R$
contains an infinite idempotent.

\begin{proposition}\label{partial converse}
If $R$ is a $1$-simple ring containing an infinite idempotent,
then it is purely infinite simple.
\end{proposition}
\begin{proof}
Let $y\in R$ be a nonzero element, and let $e\in R$ be the
infinite idempotent. By hypothesis, there exists $z,t\in R$ such
that $e=zyt$. Without loss of generality we can assume $z=ez$ and
$t=te$. Set $f=ytz$. Then, $f^2=(ytz)(ytz)=yt(zyt)z=ytez=ytz=f$,
so that it is an idempotent. Clearly, $f\in yR$, and since
$f=(yt)z$ and $e=z(yt)$, we have that $e\sim f$, whence $f$ is an
infinite idempotent, as desired.
\end{proof}

On one side, \cite[Theorem 2.2, Theorem 1.2]{LZ} imply that a
$1$-simple C*-algebra contains a nontrivial idempotent. Hence, in
the case of infinite dimensional C*-algebras, purely infinite
simple is equivalent to $1$-simple. On the other side, it is not
clear whether a $1$-simple ring has nonzero idempotents, whence
the whole equivalence remains unsolved.

\section{Algebraic Zhang's Dichotomy}

In this section we will show that an analog of Zhang's Dichotomy
for purely infinite simple C*-algebras \cite[Theorem 1.2]{Zhang}
holds in the purely algebraic context.

In order to state the results, we need to recall some definitions.
Recall that a ring $R$ is said to be exchange if for every element
$a\in R$ there exists and idempotent $e\in R$ and elements $r,s\in
R$ such that $e=ra=a+s-sa$ \cite{ext}. This definition reduces to
the Goodearl-Nicholson characterization of exchange rings in case
$R$ is a unital ring: a unital ring $R$ is said to be exchange if
for every element $a\in R$ there exists and idempotent $e\in aR$
such that $(1-a)\in(1-e)R$. Next definitions are borrowed from
\cite{AP}. Given a semiprime ring $R$, we say that a double
centralizer for $R$ is a pair $(f,g)$ such that $f:R\rightarrow R$
is a right module morphism, $g:R\rightarrow R$ is a left module
morphism, satisfying $g(x)y=xf(y)$ for all $x,y\in R$. Notice that
for any element $a\in R$, the pair $(f_a, g_a)$, where the maps
are left/right multiplication by $a$ respectively, is a double
centralizer. The set of double centralizers over $R$, endowed with
the componentwise addition and the product defined by the rule
$(f_1, g_1)\cdot (f_2,g_2)=(f_1\cdot f_2,g_2\cdot g_1)$, has
structure of ring with unit $(\mbox{Id}, \mbox{Id})$, and it is
called the ring of multipliers of $R$, denoted $\mathcal{M}(R)$.
Notice that $R$ is an ideal of $\mathcal{M}(R)$ through the
identification of $a\in R$ with $(f_a, g_a)\in \mathcal{M}(R)$;
moreover, $\mathcal{M}(R)$ coincides with $R$ whenever $R$ is a
unital ring. A net $(x_\lambda )_{\lambda \in
\Lambda}\subset\mathcal{M}(R)$ converges in the strict topology to
$x\in \mathcal{M}(R)$ if for every $a\in R$ there exists $\lambda
_0$ such that $(x_\lambda -x)a=a(x_\lambda -x)=0$ for $\lambda\geq
\lambda_0$. We say that a net $\{ a_i\}\subset R$ is an
approximate unit for $R$ provided that it converges to $1$ in the
strict topology. An approximate unit consisting on idempotents is
called a local unit. We can assume that an approximate (local)
unit is increasing \cite[Lemma 1.5]{AP}. A ring with an
approximate unit is called $s$-unital. A $s$-unital ring with a
countable approximate unit is called $\sigma$-unital. A ring has a
countable unit if it is $\sigma$-unital and has a local unit. This
is equivalent to the fact that there exists an increasing sequence
of idempotents $\{ e_n\}_{n\in \N}$ such that $R=\bigcup_{n\in
\N}e_nRe_n$ \cite[p. 3366]{AP}.

\begin{theorem}\label{spiexch} {\rm {(\cite[Theorem
1.1]{exchange})}} Every purely infinite simple ring is an exchange
ring.
\end{theorem}

We thank P. Ara for the proof of the following result.

\begin{lemma}\label{s-unit&exchange}
Every $s$-unital exchange ring is a ring with local units.
\end{lemma}
\begin{proof}
Given a finite number of elements $x_1,\dots ,x_n\in R$ we must
find an idempotent $h\in R$ such that $x_i\in hRh$ for all $i$.
Since $R$ is $s$-unital, there is $y\in R$ such that $x_iy=x_i$
for all $i$.

Let us work in $R^1=R\oplus \mathbb Z$, the unitization of $R$. By
the exchange property of $R$, there is $e=e^2\in R$ such that
$e\in yR $ and $1-e\in (1-y)R^1$. Choose $t\in R$ such that
$1-e=(1-y)(1-t)$. We then have
$$x_i(1-e)=x_i(1-y)(1-t)=0.$$
Now there is $z\in R$ such that $zx_i=x_i$ for all $i$ and $ze=e$.
Since the exchange property is left-right symmetric, there is an
idempotent $g\in R$ such that $(1-g)x_i=0$ for all $i$ and
$(1-g)e=0$. Now take $h=e+g-eg$. Then $h$ is an idempotent in $R$
and $x_i\in hRh$ for all $i$, as desired.
\end{proof}

\begin{corollary}\label{sigmaunit}
Every $\sigma$-unital exchange ring is a ring with countable unit.
\end{corollary}

The next result fills the gap to get the desired dichotomy. In
order to prove it, we need to recall a few things of K-Theory. For
a ring $R$, we denote by $V(R)$ the abelian monoid of equivalence
classes of idempotents in $M_{\infty }(R)$ under the relation
$\sim$ defined above, with the operation $[p]+[q]=[p\oplus q]$. We
consider this monoid endowed with the algebraic pre-ordering,
denoted by $\leq$, that corresponds to the ordering induced by the
relation $\lesssim$; in particular $<$ corresponds to the relation
$\prec$. Given a ring $R$, it is easy to see that $V(R)$ is
conical, and if $R$ is simple, then so is $V(R)$. If $R$ is purely
infinite simple (non necessarily unital), then the argument in the
proof of \cite[Proposition 2.1]{agp} implies that $V(R)^*$ is a
group. Hence, for every $e,f\in R$ nonzero idempotents in a purely
infinite simple ring, we have $[e]<[f]$, and thus $e\prec f$.

\begin{lemma}\label{equno}
Let $R$ be a $\sigma$-unital, non-unital, purely infinite simple
ring. For any sequence of nonzero orthogonal idempotents $\{
p_n\}_{n\geq 1}$ such that $\sum\limits_{i=1}^np_i\rightarrow P\in
\mathcal{M}(R)$ in the strict topology, $P\sim 1\in
\mathcal{M}(R)$.
\end{lemma}
\begin{proof}
By Theorem \ref{spiexch} and Corollary \ref{sigmaunit}, $R$ has a
countable unit. Let $\{ e_n\}_{n\geq 1}$ be an increasing
countable unit in $R$. Since $R$ is purely infinite simple,

$$e_1\prec p_1+p_2\prec e_3\prec p_1+p_2+p_3+p_4\prec\ldots$$

It means that there exists an idempotent $h_1\in R$ such that
$h_1\sim e_1$ and $h_1<p_1+p_2$. Hence, $p_1+p_2-h_1\prec
e_3-e_1$, and thus there exists an idempotent $g'\in R$ with
$g'\sim p_1+p_2-h_1$ and $g'<e_3-e_1$. Defining $g_2=e_1\oplus
g'\in R$, we have $e_1<e_1\oplus g'=g_2<e_1+e_3-e_1=e_3$ $g_2\sim
h_1+p_1+p_2-h_1=p_1+p_2$. By recurrence on this argument, we get
two sequences of idempotents $\{ g_{2j}\}_{j\in \N }$ and
$\{h_{2j+1}\}_{j\in \N}$ such that, for each $j\in \N$,
$e_{2j-1}<g_{2j}<e_{2j+1}$, $g_{2j}\sim p_1+\cdots +p_{2j}$, with
$p_1+\cdots +p_{2j}<h_{2j+1}<p_1+\cdots +p_{2(j+1)}$, and
$h_{2j+1}\sim e_{2j+1}$. So we have:

$$\begin{array}{ccccccccc}
h_1 & < & p_1+p_2 & < & h_3 & < & p_1+p_2+p_3+p_4 & < & \cdots \\
\wr &   & \wr     &   & \wr &   & \wr             &   &\\
e_1 & < & g_2 & < & e_3 & < & g_4 & < & \cdots
\end{array}
$$
For each $n\in \N$, define
\begin{center}
\begin{tabular}{l}
$g_n=\left\{
\begin{array}{ll}
0, & n=0; \\
e_n ,& n  \, \mbox{ odd};\\
g_n ,& n \, \mbox{ even}.
\end{array} \right.$\\
\\
$h_n=\left\{ \begin{array}{ll}
0, & n=0; \\
h_n ,& n  \, \mbox{ odd};\\
p_1+\cdots +p_n ,& n \, \mbox{ even}.
\end{array} \right.$
\end{tabular}
\end{center}
Then, we have two ascending sequences of idempotents, $\{
g_n\}_{n\in \N}$ and $\{ h_n\}_{n\in \N}$, such that $g_n\sim h_n$
for each $n\in \N$. Notice that $h_{2n}=\sum\limits_{i=1}^{2n}p_i$
in $\mathcal{M}(R)$. Also notice that, given any $a\in R$, there
exists $n\in \N$ such that, for any $m\geq n$, $h_{2m}a=Pa$. Since
$$h_{2m+2}=(h_{2m+2}-h_{2m+1})+(h_{2m+1}-h_{2m})+h_{2m},$$ defining
$ \widetilde{p}=(h_{2m+2}-h_{2m+1})$ and $
\widehat{p}=(h_{2m+1}-h_{2m})$, we have $$
\widetilde{p}a+\widehat{p}a+h_{2m}a=h_{2m+2}a=Pa=h_{2m}a.$$ Thus,
$\widetilde{p}a+\widehat{p}a=0$, and since $\widetilde{p}\perp
\widehat{p}$, we have $\widetilde{p}a=\widehat{p}a=0$. Thus, for
any $m\geq n$, $h_{2m+1}a=\widehat{p}a+h_{2m}a=Pa$. Hence,
$h_{n}\rightarrow P$ in $\mathcal{M}(R)$. Similarly we get
$g_n\rightarrow 1$ in $\mathcal{M}(R)$.

Since $R$ is purely infinite simple, $V(R)^*$ is a group
\cite[Proposition 2.1]{agp}. So, for $i \in \N$, since
$h_i+(h_{i+1}-h_i)=h_{i+1}\sim g_{i+1}=g_i+(g_{i+1}-g_i)$, we have
$h_{i+1}-h_i\sim g_{i+1}-g_i$. Thus, there exist $x_i\in
(g_{i+1}-g_i)R(h_{i+1}-h_i),\, y_i \in
(h_{i+1}-h_i)R(g_{i+1}-g_i)$ with $x_iy_i=g_{i+1}-g_i$,
$y_ix_i=h_{i+1}-h_i$. As $\{ \sum _{i=0}^n (g_{i+1}-g_i)\}_{n\in
\N}\rightarrow 1$ and $\{ \sum _{i=0}^n (h_{i+1}-h_i)\}_{n\in
\N}\rightarrow P$, by \cite[Lemma 1.7]{AP}, $\{ \sum _{i=0}^n
x_i\}_{n\in \N}\rightarrow x$ and $\{\sum _{i=0}^n y_i\}_{n\in
\N}\rightarrow y$, for some $x\in \mathcal{M}(R)P$ and $y\in
P\mathcal{M}(R)$. By \cite [Lemma 1.3]{AP}, $xy=1$ and $yx=P$.
Hence, $P\sim 1$ in ${\mathcal{M}(R)}$.
\end{proof}

Finally, we get the main result in this paper.

\begin{theorem}\label{dicotomia}
Let $R$ be a $\sigma$-unital, non-unital, purely infinite simple
ring. Then:
\begin{enumerate}
\item $R\cong M_\infty (R)$; \item For every nonzero idempotent
$q\in R$, we have $R\cong M_\infty (qRq)$.
\end{enumerate}
\end{theorem}
\begin{proof}
By Theorem \ref{spiexch} and Corollary \ref{sigmaunit}, $R$ has a
countable unit. Let $\{ e_n\}_{n\geq 1}$ be an increasing
countable unit in $R$. Fix a nonzero idempotent $q\in R$. We
define a sequence of idempotents by recurrence, as follows:
\begin{center}
\begin{tabular}{l}
$q_0=0;$\\ $q_n=e_n-e_{n-1}\mbox{, } n\in \N \quad (e_0=0)$
\end{tabular}
\end{center}

Since $R$ is purely infinite simple, $q_n$ is an infinite
idempotent for any $n\in \N$. Moreover, $q\lesssim q_n$. Hence,
for each $n\in \N$, there exists an idempotent $p_n\in R$ such
that $p_n\leq q_n$ and $p_n\sim q$.

By construction, $e_n=\sum_{i=0}^nq_i$, and $\{ e_n\}_{n\in \N}=\{
\sum_{i=0}^nq_i\}_{n\in \N}$ converges to $1\in \mathcal{M}(R)$ in
the strict topology of $R$; in particular, it is a Cauchy
sequence. Since $R$ is simple, it is semiprime, and
$$(\sum_{i=0}^np_i-\sum_{i=0}^mp_i)=\sum _m^np_i\leq \sum_m^nq_i$$
implies that $(\sum_m^np_i)_{n\in \N}$ is also a Cauchy sequence.
By \cite[Proposition 1.6]{AP}, $ \mathcal{M}(R)$ is complete, so
that $(\sum_m^np_i)_{n\in \N}$ converges to some $P\in
\mathcal{M}(R)$.

Clearly, $\{q_n\} _{n\in \N}$ is a family of orthogonal
idempotents. Then, by \cite[Lemma 1.3]{AP}
$$P^2=(\lim_n\sum_{i=1}^np_i)(\lim_n\sum_{j=1}^np_j)=
\lim_n(\sum_{i=1}^np_i)(\sum_{j=1}^np_j)=
\lim_n\sum_{i=1}^np_i=P,$$

\noindent whence $P$ is an idempotent of $\mathcal{M}(R)$. By
Lemma \ref{equno}, $P\sim 1\in \mathcal{M}(R)$. In particular,
there exist $u\in P\mathcal{M}(R)$ and $v\in \mathcal{M}(R)P$ such
that $uv=P, vu=1$. Notice that, since $R$ is non-unital, $P
\not\in R$.

Thus, we can define two ring morphisms, $\rho: R\rightarrow PRP$
by the rule $\rho (r)=urv$, and $\psi: PRP\rightarrow R$ by the
rule $\psi (r)=vru$. Clearly they are mutually inverses, so that,
\begin{equation}\label{isom1}
R\cong PRP.
\end{equation}

Define $t_n=\sum _{i=1}^n p_i$. Since $\{ t_n\}_{n\in \N}$
converges to $P\in \mathcal{M}(R)$, we have $PRP=\bigcup_{n\in
\N}t_nRt_n$. Then, $t_{n+1}-t_n=\sum _{i=1}^{n+1} p_i-\sum
_{i=1}^n p_i=p_{n+1}\sim q$, and since
$t_n=(t_n-t_{n-1})\oplus(t_{n-1}-t_{n-2}) \oplus \ldots \oplus
(t_1-t_0)\sim nq$, we get

$$t_nRt_n\cong End_R(t_nR)\cong End_R(n(qR))\cong M_n(qRq).$$

Under this identification, $t_nRt_n\hookrightarrow
t_{n+1}Rt_{n+1}$ is the map
\begin{center}
\begin{tabular}{ccc}
$M_n(qRq)$ & $\longrightarrow$ & $M_{n+1}(qRq)$\\
 $a$ & $\longmapsto$ & $diag(a,0)$
\end{tabular}
\end{center}
so that
\begin{equation}\label{isom2}
PRP=\bigcup_{n\in \N}t_nRt_n\cong\bigcup_{n\in
\N}M_n(qRq)=M_\infty(qRq).
\end{equation}

Finally, if $q\in R$ is a nonzero idempotent, $qRq$ is a unital,
purely infinite simple ring. Then, (\ref{isom1}) and (\ref{isom2})
imply $R\cong PRP\cong M_\infty (qRq)$. Hence, $M_\infty (R)\cong
M_\infty (M_\infty (qRq))\cong M_\infty (qRq)\cong R$, as desired.
\end{proof}

Then, we get the corresponding Dichotomy result, analog to
\cite[Theorem 1.2(i)]{Zhang}. We say that a (non-unital) ring $R$
is stable if there exists a ring $S$ such that $R\cong M_{\infty
}(S)$.

\begin{corollary}
Let $R$ be a $\sigma$-unital purely infinite simple ring. Then it
is either unital or stable.
\end{corollary}

\begin{remark}\label{dicolimit}
{\rm Notice that we cannot guarantee that a non-unital, purely
infinite simple ring has $s$-unit. For example, given a field $K$,
consider, for $n\geq 2$, the Leavitt algebra
$$R=K\langle x_1,\ldots ,x_n,y_1,\ldots ,y_n  \,| \,
x_iy_j=\delta _{ij}, \sum_{i=1}^ny_ix_i=1\rangle.$$ This is a
purely infinite simple ring (see \cite{agp}), so that any right
ideal of $R$ is a non-unital purely infinite simple ring. Then, it
is easy to see that the right ideal $L=y_1R$ is a non-unital,
purely infinite simple ring with no $s$-unit.}
\end{remark}

\section*{Acknowledgments}

We thank P. Ara for valuable discussions during the preparation of
this note.

\end{document}